\documentclass[12pt]{amsart}
\usepackage{geometry}                % See geometry.pdf to learn the layout options. There are lots.
\usepackage{color}\usepackage{pdfpages}
\usepackage{fullpage}
\usepackage{wrapfig}%
\usepackage{caption}

\newcommand{\Order}[1]{\ensuremath{\mathcal{O}(#1)}}    % big O notation

\newcommand{\appropto}{\mathrel{\vcenter{
  \offinterlineskip\halign{\hfil$##$\cr
    \propto\cr\noalign{\kern2pt}\sim\cr\noalign{\kern-2pt}}}}}
    
\begin{document}

\title{Segmental Refinement: A Multigrid Technique for Data Locality}

\author{Mark F. Adams}
%\email{\tt mfadams@lbl.gov}
\thanks{Scaleable Solvers Group, Lawrence Berkeley National Laboratory, Berkeley, CA}
\author{Jed Brown}
\thanks{
%Laboratory for Advanced Numerical Simulation Mathematics and Computer Science Division, Argonne National Laboratory, Department of Computer Science University of Colorado Boulder
Mathematics and Computer Science Division, Argonne National Laboratory
}
\author{Matt Knepley}
\thanks{Computation Institute, University of Chicago}
\author{Ravi Samtaney}
\thanks{King Abdullah University of Science and Technology}
%\date{}                                           % Activate to display a given date or no date


\date{}

%
% Tau correction,  Achi uses even 2J(k) >= s(M-k+1), s == min(smoothness order, discretization order) (2), "segmental refinement" algorithms
%

\begin{abstract} 

We investigate a domain decomposed multigrid technique, segmental refinement, for solving general nonlinear elliptic boundary value problems. 
Brandt and Diskin first proposed this method in 1994;  we continue this work by analytically and experimentally investigating its complexity.
We confirm that communication of traditional parallel multigrid can be eliminated on fine grids with modest amounts of extra work and storage while maintaining the asymptotic exactness of full multigrid, although we observe a dependence on an additional parameter not considered in the original analysis.
We present a communication complexity analysis that quantifies the communication costs ameliorated by segmental refinement and report performance results with up to 64K cores of a Cray XC30.
 
\end{abstract}

\maketitle

\section{Introduction}

Full multigrid (FMG) is a provably asymptotically exact, non-iterative algebraic equation solver for discretized elliptic partial differential equations (PDEs) with work complexity of about five residual calculations, or what is known as textbook multigrid efficiency, for the constant coefficient Laplacian \cite{Bank-81}.
While textbook multigrid efficiency is only provable for a few classes of elliptic problems, it has been observed experimentally in many more problems \cite{Thomas2001853,UTrottenberg_CWOosterlee_ASchueller_2000a,Adams-10a}, and is applicable to general nonlinear elliptic equations.
Multigrid methods are widely used in practice; they are important methods to adapt to emerging architectures. % that we are encountering today and anticipate in the future.

Memory movement, in both intra-node and inter-node communication, and global data dependencies are the primary drivers of costs, in power and time, for PDE simulations on current, and anticipated future, computer architectures.
Memory movement pressures are not new and have been accumulating for decades, but the recent prominence of energy costs in powering memory and moving data is exacerbating this problem.
Segmental refinement addresses the challenges posed by the deep memory hierarchies of modern architectures at a fundamental, algorithmic level by exploiting the local nature of multigrid processes and a tolerance for finite algebraic error, which nonetheless vanishes asymptotically.
A segmental refinement data model or method explicitly decouples subdomain processing, at some level of the memory and multigrid hierarchy, which improves data locality, amortizes latency costs, and reduces data dependencies.

Brandt proposed segmental refinement in the 1970s \cite{Brandt-77} \S 7.5;  \cite{Brandt-2011} \S 8.7; \cite{NDinar_1979a} as a low memory complexity technique for FMG that does not store the entire solution in memory at any one time.
Brandt and Diskin recognized that segmental refinement has attractive properties for distributed memory computing \cite{ABrandt_BDiskin_1994a}; it is inherently asynchronous and highly parallel, with no interprocess communication on the finest grids, and it requires only modest amounts of extra storage and work in buffer cells.
This paper continues the development of segmental refinement by quantifying its complexity both experimentally and analytically.
We present the first published multilevel numerical results, and report preliminary performance results, of a cell centered segmental refinement data model.

\section{Differential and discretized problems}
\label{sec:defs}

We consider general nonlinear elliptic problems in an open domain $\Omega$ with boundary $\partial \Omega$ of the form
\begin{equation} 
\label{eq:eq}
Lu\left(x\right) = f\left(x\right)  \qquad \left(x \in \Omega\right),
\end{equation}
where $f$ is a known function, $u$ is unknown, and $L$ is a uniformly elliptic operator.
While the methods described herein are generally applicable, we restrict ourselves to the 3D Poisson operator: $x=(x_1,x_2,x_3)$, $L=\left( \frac{\partial^2 u}{\partial x_1^2}
      + \frac{\partial^2 u}{\partial x_2^2}
      + \frac{\partial^2 u}{\partial x_3^2} \right)$.
In addition to this interior equation, a suitable boundary condition on $\partial \Omega$ is assumed; we assume $u(x) = 0$; $x \in \partial  \Omega$.

The discretization of equation (\ref{eq:eq}) can be fairly general but we restrict ourselves to a cell-centered finite difference method on isotropic Cartesian grids and rectangular domains.
For the grid $\Omega_h$ with mesh spacing $h$ covering the domain $\Omega$, the equation can be written as
\begin{equation} 
\label{eq:deq}
L_hu_h(i) = f_h(i)  \qquad \left(i \in \Omega_h\right),
\end{equation}
where $i = \left( x - \frac{h}{2} \right) / h$ is an integer vector, $x = ih + \frac{h}{2}$ is a cell center, and the boundary $\partial \Omega$ lines up with the cell edges.
In 3D $i = \left(i_1,i_2,i_3\right)$ is an index for a cell in grid $\Omega_h$. 
The indexing in equation (\ref{eq:deq}) is dropped and field variables (e.g., $u_h$) are vectors of scalars.

Our grids $\Omega_h$ and subsequent subdomains are isotropic and can for the most part be expressed as tensor products of 1D grids.
The lengths of $\Omega_h$, in each dimension, is an integer vector; we simplify the presentation by using the integer $N$, because we use cubical subdomains.
Multigrid utilizes an accurate and inexpensive  solver on $\Omega_0$ and a sequence of grids $\Omega_0$, $\Omega_1$, $\Omega_2$,..., $\Omega_M$, where $\Omega_k \equiv \Omega_{h_k}$, $h_k = h_{k-1}/2$, $h_M=h$, $N_k = 2N_{k-1}$, $N_M = N$.

\section{Multigrid background}

The antecedents of modern multigrid go back at least to Southwell in the 1930s \cite{RVSouthwell_1940a}, Fedorenko in the early 1960s \cite{RPFedorenko_1961a}, and others \cite{UTrottenberg_CWOosterlee_ASchueller_2000a}.
Brandt developed multigrid's modern form in the 1970s, an asymptotically exact solver with work complexity of a few residual calculations -- what is known as textbook multigrid efficiency.
He applied multigrid to complex domains, variable coefficients, and nonlinear problems \cite{ABrandt_1973a}.
A substantial body of literature, both theoretical and experimental, demonstrates the efficacy of multigrid \cite{UTrottenberg_CWOosterlee_ASchueller_2000a,Brandt-2011}.
Full Approximation Scheme (or Storage, FAS) multigrid has also been demonstrated to be an effective nonlinear solver, with costs  similar to those of a linearized multigrid solve (e.g., \cite{UTrottenberg_CWOosterlee_ASchueller_2000a} \S 5.3.3, \cite{Adams-10a}).

\subsection{Multigrid algorithm}

Multigrid starts with the observation that errors that are poorly resolved with local processes can often be resolved with local processes on a lower resolution discretization.
This lower resolution problem is known as a coarse grid.
Multigrid, by definition, applies this process recursively until the problem size is small enough to be solved inexpensively and exactly.
The coarse grid space is represented algebraically by the columns of the prolongation operator $I^h_{H}$ or $I^k_{k-1}$, where $h$ is the fine grid mesh spacing and $H$ is the coarse grid mesh spacing.
Residuals are mapped from the fine grid to the coarse grid with the restriction operator $I^H_{h}$. 
The coarse grid operator can be formed in one of two ways (with some exceptions), either algebraically to form Galerkin (or variational) coarse grids, $L_{H} = I^H_{h}L_{h}I^h_{H}$, or by creating a new operator on each coarse grid if an explicit coarse grid with boundary conditions is available.

Correction Scheme (CS) multigrid, where coarse grids compute corrections to the solution, is appropriate for linear problems, but FAS multigrid is more natural for segmental refinement.
FAS is derived by writing the coarse grid residual equation for equation (\ref{eq:deq}) as
\begin{equation} 
r_{H} = L_{H} (u_{H} ) - L_{H} ({\hat u}_{H} ) = L_{H} ({\hat u}_{H} +e_{H} ) - L_{H} ({\hat u}_{H} ),
\label{eq:cresid}
\end{equation}
where $u_H$ is the exact solution, ${\hat u}_{H}$ is an approximation to $I^H_h{u}_h$ (which is the full solution represented on the coarse grid), and $e$ is the error.
With an approximate solution on the fine grid $\tilde u_h$, the coarse grid equation can be written as
\begin{equation} 
L_{H}\left (I^{H}_h {\tilde u}_h+ e_{H}\right) = L_{H} \left(I^{H}_h{\tilde u}_h\right) + I^{H}_h\left( f_h - L_h {\tilde u}_h \right) = f_H =  I^{H}_h\left( f_h\right) + \tau^{H}_h, 
\label{eq:cresid2}
\end{equation}
and is solved approximately; $\tau^{H}_h$ is the tau correction, which represents a correction to the coarse grid from the fine grid.
After $I^{H}_h {\tilde u}_h$ is subtracted from the $I^{H}_h {\tilde u}_h+ e_{H}$ term the correction is applied to the fine grid with the standard prolongation process.
Figure \ref{fig:mgv} shows an FAS multigrid $V(\nu 1,\nu 2)$--cycle algorithm with nonlinear local process or smoother $u \leftarrow S(L,u,f)$.
%\vskip .2in
\begin{figure}[h]
\vbox{ \raggedright 
$\phantom{}u=$ {\bf  function} $FASMGV(L_k,u_k,f_k)$ \\ 
$\phantom{MM}${\bf if $k > 0$} \\
$\phantom{MMMM}u_k \leftarrow S^{\nu 1}(L_k,u_k,f_k)$ \\
$\phantom{MMMM}r_k \leftarrow f_k - L_ku_k$ \\ 
$\phantom{MMMM}u_{k-1}\leftarrow {\hat I}^{k-1}_k(u_k)$ \\ 
$\phantom{MMMM}r_{k-1}\leftarrow I^{k-1}_k(r_k)$ \\ 
$\phantom{MMMM}t_{k-1}\leftarrow u_{k-1}$ \\ 
$\phantom{MMMM}w_{k-1}\leftarrow FASMGV(L_{k-1},u_{k-1},r_{k-1}+L_{k-1}u_{k-1})$  \\
$\phantom{MMMM}u_k \leftarrow u_k + I^k_{k-1}(w_{k-1} - t_{k-1})$ \\ 
$\phantom{MMMM}u_k \leftarrow S^{\nu 2}(L_k,u_k,f_k)$ \\ 
$\phantom{MM}${\bf else}\\ 
$\phantom{MMMM}u_k \leftarrow L_k^{-1}f_k$ \\ 
$\phantom{MM}${\bf return} $u_k$}
\caption{FAS multigrid {\it $V$-$cycle$} }
\label{fig:mgv}
\end{figure}

A lower order restriction operator,  ${\hat I}^{H}_h$, can be used to restrict solution values if a higher order ${I}^{H}_h$ is used for the residual, because this approximate coarse grid solution is subtracted from the update to produce an increment and is only needed for the  nonlinearity of the operator (e.g., ${\hat I}^{H}_h=0$ recovers CS multigrid).
 
\subsection{Full Multigrid Algorithm}
\label{ssec:fmg}
An effective V--cycle reduces the error by a constant fraction and is thus an iterative method, but it can be used to build a non-iterative, asymptotically exact solver that reduces the algebraic error to the order of the discretization error.
FMG starts on the coarsest grid where an inexpensive accurate solve is available, prolongates the solution to the next finest level, applies a V--cycle, and continues until a desired resolution is reached.
Figure \ref{fig:fmg} shows the full multigrid algorithm, with $M$ coarse grids and $\alpha$ steps of the smoother before each V--cycle, in an F($\alpha$,$\nu 1$,$\nu 2$) cycle.
%\vskip .2in
\begin{figure}[h]
\vbox{ \raggedright 
$\phantom{}u=${\bf  function} $FMG$ \\ 
$\phantom{MM}u_{0} \leftarrow 0$ \\
$\phantom{MM}u_{0} \leftarrow FASMGV\left(L_{0}, u_{0}, f_{0}\right)$  \\
$\phantom{MM}${\bf for k=1:M} \\
$\phantom{MMMM}{u}_{k} \leftarrow \Pi^{k}_{k-1} u_{k-1}$   \\
$\phantom{MMMM}{u}_{k} \leftarrow S^{\alpha}(L_k,u_k,f_k)$  \\
$\phantom{MMMM}{u}_{k} \leftarrow FASMGV\left(L_k,u_k,f_k\right)$ \\
$\phantom{MM}${\bf return} $u_0$}
\caption{ Full multigrid}
\label{fig:fmg}
\end{figure}
A higher order interpolator between the level solves, $\Pi^h_{H}$, is useful for optimal efficiency of FMG and is required if ${I}^{H}_h$ is not of sufficient order (e.g., $\Pi^h_{H}$ must be at least linear, for  cell-centered $2^{nd}$-order accurate discretizations, whereas ${I}^{H}_h$ can be constant).

One can analyze FMG with an induction hypothesis that the ratio $r$ of the algebraic error to the discretization error is below some value and assume that the discretization error is of the form $Ch^{p}$, where $p$ is the order of accuracy of the discretization.
Further, assume that the solver on each level (e.g., one V--cycle) reduces the error by some factor $\Gamma$, which can be proven or measured experimentally, to derive the relationship between $\Gamma$ and $r$: $\Gamma = \frac{r}{\left(4r+3\right)}$, with $p=2$ and a refinement ratio of two.
One must use a sufficiently powerful solver such that $\Gamma < 0.25$.
For instance, Adams et. al. use FMG for compressible resistive magnetohydrodynamics problems; two V--cycles were required to achieve sufficient error reduction $\Gamma$ \cite{Adams-10a}.

\subsection{Conventional distributed memory multigrid}
\label{sec:pmg}

Domain decomposition is a natural technique for distributed memory processing of many classes of discretized PDEs, where each subdomain is placed on a processor or memory partition and the semantics of the serial algorithm are replicated.
This process starts by decomposing $\Omega_h$ into $P$ disjoint grids ${}^p\Omega_h$ such that $\Omega_h = \bigcup_{p=1}^P {}^p\Omega_h$.
We use a rectangular array of processes of size $\left(P_1,P_2,P_3\right)$ and thus $P=P_1P_2P_3$.
We implement boundary conditions in equation (\ref{eq:deq}) with ghost cells; $\Omega_h$ is enlarged by one cell in all directions to form $\Omega_h^{+1}$.
Boundary ghost cell values are set with appropriate (linear) interpolation of interior values before each operator application.

We define the number of cells on each side of a (cube) subdomain ${}^p\Omega_h$ as the integer ${}^pN_k$ on level $k$, again using integers for simplicity.
The total number of cells in our problems is thus $n=P_1P_2P_3{}^pN^3_M$, where ${}^pN_M$ is the number of cells in each dimension on the fine grid.
A conventional distributed memory full multigrid algorithm starts with a small coarse grid on a small number of processes (e.g., one process).
The coarse grid is refined and split into equally sized patches, which populate more processes.
This process continues until all processes are used, forming an octree in 3D.
We continue with simple refinement once all processes are used; however, more complex distributed memory models are common \cite{HPGMGv1}.

%\pagebreak
\section{Segmental refinement}
\label{sec:sr}

%\begin{figure}[h]
\begin{wrapfigure}{R}{0.5\textwidth}
\begin{center}
\includegraphics[width=80mm]{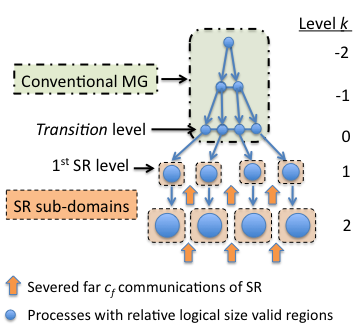}
\caption{1D SR data model}
\label{fig:1dsr-a}
\end{center}
%\end{figure}
\end{wrapfigure}

This section describes a cell-centered segmental refinement (SR) data model or method.
% that is a straightforward adaptation of the vertex-centered model of Brandt and Diskin.
Segmental refinement begins with a conventional distributed memory FAS-FMG method, which is used as a ``coarse" grid solver. % and as the basis for the SR grids.
The finest level of this solver is the ``transition" level and is given a grid index $k=0$; coarser grids have negative indices.
The subsequent $K$ fine grids are defined as SR grids.
Figure \ref{fig:1dsr-a} shows a 1D example, with two SR levels and four processes.
%This is a simple SR data model because the transition level and SR levels use the same parallel decomposition and no distributed memory parallelism is used within the SR subdomains.

Non-ghost cells are defined as genuine cells and the {\it genuine region} is defined as ${}^p\Omega_k^V \equiv ~^p\Omega_k$.
SR adds buffer cells by growing each local subdomain grid by $2J_k$ cells in each dimension.
Following Brandt and Diskin, we define the length, in each dimension,  of the SR buffer region $J_k$ to be
\vspace{.5em}
\begin{equation} 
J_k = J(k) = 2 \cdot \Bigl\lfloor \frac{A + B \cdot \left(K-k\right)}{2}\Bigr\rfloor,
\label{eq:nb}
\end{equation}
where $A$ is a constant term and $B$ is a linearly increasing term on coarser grids.
$J_k$ is constrained to even integers to simplify restriction.
The union of the genuine cells and the $J_k$ buffer cells defines the {\it compute region} ${}^p\Omega_k^C \equiv {}^p\Omega_k^{V+J_k} \cap \Omega_k$, that is, the genuine region grown by $J_kh_k$ in all directions and clipped by the domain.
The length of the compute region is generally ${}^pN^C_k = {}^pN^V_k + 2J_k$, where ${}^pN^V_k$ is the length of the compute region in grid $k$. % away from domain boundaries.

%\begin{figure}[h]
\begin{wrapfigure}{R}{0.5\textwidth}
%\begin{center}
\includegraphics[width=80mm]{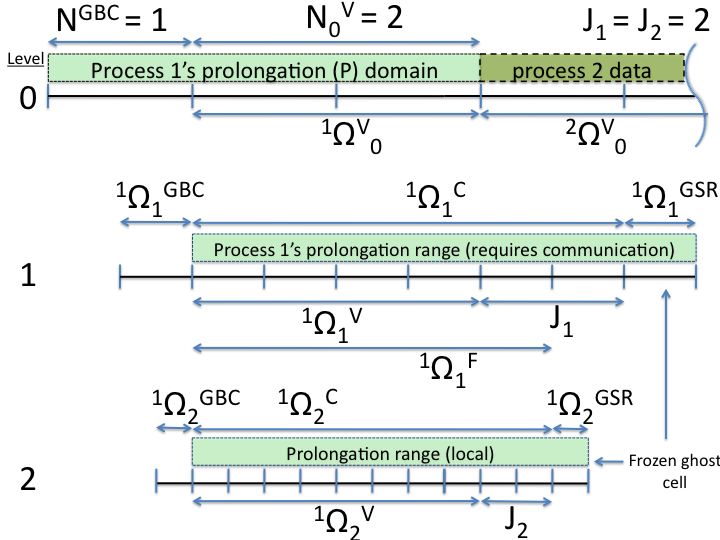} 
\caption{1D 2 process SR, with the number of boundary condition cells $N^{GBC} $ }
\label{fig:SR-1D}
%\end{center}
%\end{figure}
\end{wrapfigure}c

We define {\it process ghost cells} by ${}^p\Omega_k^G \equiv {}^p\Omega_k^{C+1} \setminus  {}^p\Omega_k^{C}$ and subdivide ${}^p\Omega_k^G$ into two sets: ${}^p\Omega_k^{GBC} \equiv {}^p\Omega_k^{G} \setminus \Omega_k$ and ${}^p\Omega_k^{GSR} \equiv {}^p\Omega_k^{G} \cap \Omega_k$.
The ${}^p\Omega_k^{GBC}$ cell values are computed with the conventional boundary condition algorithm.
${}^p\Omega_k^{GSR}$ cells are set during the $I_0^1$ prolongation process and are ``frozen", in that they are not updated with the neighbor exchanges, during the rest of the multigrid process.
This ``freezing" is a consequence of the elided communication of SR; ${}^p\Omega_k^{GSR}$ cells are set with prolongation only. 
Define the {\it support of the compute region}, on grid $k$ , of grid $k+1$ as ${}^p\Omega_k^F \equiv {}^p\Omega_{k+1}^C$; this is the region updated with the simple averaging restriction operator.
Figure \ref{fig:SR-1D} shows a 1D example at the edge of the domain with two processes and two SR levels with the range of prolongation for one process. % and ghost cells.

The $\tau$ correction is modified to accommodate the lack of an update in the region ${}^p\Omega^C_k \setminus  {}^p\Omega^F_k$.
The range of prolongation is ${}^p\Omega^C \cup {}^p\Omega^{GSR}$.
Figures \ref{fig:fmgsr} and \ref{fig:mgvsr} show the SR FAS-FMG algorithm with annotations for the domain of each operation.
\begin{figure}[h]
\vbox{ \raggedright 
$\phantom{}u=${\bf  function} $FASFMGSR$ \\ 
$\phantom{MM}u_0 \leftarrow FMG$  \\ 
$\phantom{MM}${\bf for k=1:K} \\
$\phantom{MMMM}{u}_{k} \leftarrow \Pi^{k}_{k-1} u_{k-1}$ \qquad  \qquad  \qquad \qquad ${}^p\Omega_k^{C} \cup {}^p\Omega^{GSR}_{k}$   \\ 
$\phantom{MMMM}{u}_{k} \leftarrow S^{\alpha}(L_k,u_k,f_k)$ \qquad \qquad \qquad ${}^p\Omega^C_{k}$\\
$\phantom{MMMM}u_{k}\leftarrow FASMGVSR(L_{k},u_{k},f_{k})$ \quad   \\ 
$\phantom{MM}${\bf return} $u_K$}
\caption{FMG segmental refinement}
\label{fig:fmgsr}
\end{figure}

\begin{figure}[h]
\vbox{ \raggedright 
$\phantom{}u=${\bf  function} $FASMGVSR(L_k,u_k,r_k)$ \\ 
$\phantom{MM}u_k \leftarrow S^{\nu 1}(L_k,u_k,r_k)$, \qquad \qquad \qquad  \quad  \qquad \qquad \quad${}^p\Omega^C_{k}$ \\ 
$\phantom{MM}u_{k-1} \leftarrow {\hat I}^{k-1}_k(u_k)$ \qquad \qquad \qquad \qquad \qquad \qquad \qquad ${}^p\Omega^F_{k}$ \\ 
$\phantom{MM}t_{k-1} \leftarrow u_{k-1}$ \\
$\phantom{MM}r_{k-1} \leftarrow I^{k-1}_k(r_k - L_ku_k) + L_{k-1} u_{k-1}$, \qquad \qquad \quad ${}^p\Omega^F_{k}$ \\ 
$\phantom{MM}r_{k-1} \leftarrow L_{k-1} u_{k-1}$,  \qquad \qquad \quad \qquad \qquad \qquad \qquad ${}^p\Omega^C_{k} \setminus {}^p\Omega^F_{k}$ \\ 
$\phantom{MM}${\bf if $k = 1$} \\
$\phantom{MMMM}w_{k-1} \leftarrow FASMGV(L_{k-1},u_{k-1},r_{k-1})$ \\
$\phantom{MM}${\bf else} \\
$\phantom{MMMM}w_{k-1} \leftarrow FASMGVSR(L_{k-1},u_{k-1},r_{k-1})$ \\
$\phantom{MM}u_k \leftarrow u_k + I^k_{k-1}(w_{k-1} - t_{k-1})$, \qquad  \qquad \qquad \qquad   ${}^p\Omega^C_{k} \cup {}^p\Omega^{GSR}_{k}$ \\ 
$\phantom{MM}u_k \leftarrow S^{\nu 2}(L_k,u_k,r_k)$, \qquad \qquad \quad \qquad \quad  \qquad  \qquad  ${}^p\Omega^C_{k}$\\
$\phantom{MM}${\bf return} $u_k$}
\caption{FAS V--cycle segmental refinement}
\label{fig:mgvsr}
\end{figure}

%\pagebreak
\section{Experimental observation of parameter requirements}
\label{ssec:exp_params}

This section experimentally investigates the parameters required to maintain an acceptably accurate segmental refinement FMG solver.
There are several parameters that define the SR solver: the number of SR levels $K$ and the total number levels $M+1$; $A$ and $B$ of equation (\ref{eq:nb}); and the length of the subdomains on the transition level ${}^pN_0^V$.

\subsection{Model problem and solver}
\label{ssec:problem}
We use a multigrid refinement ratio of two, piecewise constant restriction, and linear prolongation for both the FMG and V--cycle prolongation.
The pre- and post-smoothers are $2^{nd}$-order Chebyshev polynomials and the pre V--cycle smoother is a $1^{st}$-order Chebyshev polynomial (an F(1,2,2) cycle).
The solution is prescribed as $u=\prod\limits_{i=1}^3\left( x_i^4 - R_i^2x_i^2 \right)$, for the Laplacian $Lu=f$, on a rectangular domain $$\Omega = \left\{ x_1, x_2, x_3 \ge 0, x_1 \le 2, x_2, x_3 \le 1\right\}, $$ with $R=(2,1,1)$ and a 4~x~2~x~2 process grid.
We use a homogenous Dirichlet boundary condition and a 27-point finite volume stencil that is $2^{nd}$-order accurate.

\subsection{Experiments}
We define an acceptable level of error, in the infinity norm, to be less than about $10\%$ more than the conventional solver error ($e_{conv}$); the conventional solver is $2^{nd}$-order convergent.
%We sample the high dimensional parameter design space, $K$, $A$, $B$, and ${}^pN_0^V$, to experimentally determine the manifold at which the solver error transitions from acceptable to unacceptable.
We sample the parameter space of $K$, $A$, $B$, ${}^pN_0^V$, to find the manifold where the solver error transitions from acceptable to unacceptable.
%We require that the solver maintain $2^{nd}$ order convergence rates in the max norm of the error as the grid is refined.
Table \ref{tab:paramsABC} shows the ratio ($e_r$) of the SR error ($e_{SR}$) to $e_{conv}$ ($e_r \equiv e_{SR}/e_{conv}$) with $A=2,4,6,8$ (tables), $B=0,1,2,3$ (rows), and $\log_2 {}^pN_0^V$ and $K$ (columns), and underlines the largest acceptable point in each column.
\begin{table}[h!]
\begin{minipage}[t]{.24\linewidth}
\begin{tabular}{|c||c|c|c|}\hline
& \multicolumn{3}{c|}{$\log_2 {}^pN_0^V$ ($K$)} \\\hline\hline
B &  4(6)  & 3(5) & 2(4) \\\hline
0   &  {17}    &  {7.2}  &  {2.7}  \\
1   &  {2.9}   &  {2.1}  &  {1.2}  \\
2   &  {1.5}   &  {1.2}  &  {NA}  \\
3   &  {1.2}   &  \underline{1.1}  &  {NA}  \\
\hline 
\end{tabular}
\center{(a) A=2}
\end{minipage}
\begin{minipage}[t]{.24\linewidth}
\begin{tabular}{|c||c|c|c|}\hline
  & \multicolumn{3}{c|}{$\log_2 {}^pN_0^V$ ($K$)} \\\hline\hline
B &  4(6)  & 3(5) & 2(4) \\\hline
0  &  {5.7}  & {2.6} & {1.2}  \\
1  &  {2.0}  & {1.4} & \underline{1.0}  \\
 2  & {1.3}  & \underline{1.1} & {NA}  \\
 3  & \underline{1.1}  & {1.0} & {NA}  \\
\hline
\end{tabular}
\center{(b) A=4}
\end{minipage}
\begin{minipage}[t]{.24\linewidth}
\begin{tabular}{|c||c|c|c|}\hline
& \multicolumn{3}{c|}{$\log_2 {}^pN_0^V$ ($K$)} \\\hline\hline
B &  4(6)  & 3(5) & 2(4) \\\hline
0  &  {2.8}  & {1.4} & \underline{1.0}  \\
1  &  {1.5}  & \underline{1.1} & {NA}  \\
 2  & {1.3}  & {1.0} & {NA}   \\
 3  & \underline{1.1}  & {NA} & {NA}   \\
\hline
\end{tabular}
\center{(c) A=6}
\end{minipage}
\begin{minipage}[t]{.24\linewidth}
\begin{tabular}{|c||c|c|c|}\hline
  & \multicolumn{3}{c|}{$\log_2 {}^pN_0^V$ ($K$)} \\\hline\hline
B &  4(6)  & 3(5) & 2(4) \\\hline
0  &  {1.5}  & \underline{1.1} & \underline{1.0}  \\
1  & {1.3}  & {1.0} & {NA}  \\
 2  & \underline{1.1}  & {1.0} & {NA}  \\
 3  & {1.0}  & {NA} & {NA}   \\
\hline
\end{tabular}
\center{(d) A=8}
\end{minipage}
\caption{$e_r$: ratio of SR to conventional multigrid solution error}
\label{tab:paramsABC}
\end{table}
The total number of multigrid levels can be inferred from ${}^pN_0^V$ and the process grid (i.e., $M = K + \log_2 {}^pN_0^V +  \log_2 {P_z} = K + \log_2 {}^pN_0^V +  2 $).
This data shows that $A$ and $B$ both correlate with increased accuracy, which is expected because they both increase $J$.
We observe that doubling the length of ${}^pN_0^V$ with $K$ ($\log_2 {}^pN_0^V \propto K$) and increasing $B$ with $K$ ($B \propto K$) appears to maintain an asymptotically exact solver; we base the design of a putative asymptotically exact SR data model on this observation in \S\ref{sec:new_sr_model}.
% although this data is not conclusive.
%This behavior was not anticipated from previous work on SR, but again these are the first known multilevel numerical results.
%Future work includes verifying these results, using vertex centered finite element discretizations, and higher order prolongation.% on the asymptotics of SR.

To further investigate the effect of ${}^pN_0^V$ on error we fix $A=8$, $B=0$ ($J_k=8$), and $K=5$; the relative error as a function of ${}^pN_0^V$ is shown in Table \ref{tab:C}.
\begin{table}[h!]
\begin{center}
\begin{tabular}{|c||c|c|c|c|}\hline
$\log_2 {}^pN_0^V$ &  5       & 4     & 3     & 2 \\
$N = N_K$ &  1024  & 512 & 256 & 128  \\\hline
$e_r$ &   1.02  & 1.05 & 1.13 & 1.28 \\\hline
\end{tabular}
\end{center}
\caption{$e_r({}^pN_0^V)$: $A=8, B=0, K=5$}
\label{tab:C}
\end{table}
This data shows a reduction in the error by a factor of about two with a doubling of ${}^pN_0^V$.

\subsubsection{Maximum segmental refinement buffer schedule}

The buffer length of the coarsest SR grid, $J_1$, is an important parameter because these cells require communication; these cells are the range of prolongation to the coarsest SR level and the data source for all subsequent finer grid processing.
To investigate the relationship of $J_1$ to accuracy we test with a maximum buffer schedule (MBS), where $J_1$ is a parameter and rest of the SR buffers completely support grid $k=1$: ${}^p\Omega^C_k = {}^p\Omega^F_k$ for $k < K$.
This is probably not a practical buffer schedule, because the number of buffer cells increases exponentially with refinement; the MBS removes one source of error: the lack of update of the solution and $\tau$ correction in ${}^p\Omega^C_{k} \setminus {}^p\Omega^F_{k}$.
Table \ref{tab:C_maxB} shows the error ratio as a function of ${}^pN_0^V$ with fixed $K=4$ and $J_1=4$ and the maximum buffer schedule.
\begin{table}[h!]
\begin{tabular}{|c||c|c|c|c|c|c|}\hline
$\log_2 {}^pN_0^V$ & 6       & 5     & 4      & 3     & 2     & 1\\
 $N = N_K$        & 1024 & 512 & 256  & 128 & 64   & 32 \\\hline
$e_r$  & 1.02 & 1.05 & 1.11 & 1.25 & 1.4 & 1.9 \\\hline
\end{tabular}
\caption{Effect of ${}^pN_0^V$ on error ratio with $J_1=4$, $K=4$, with MBS}
\label{tab:C_maxB}
\end{table}

This data shows slightly less degradation of the solution with increasing ${}^pN_0^V$ than that of Table \ref{tab:C}, but we observe a similar doubling of the error with each halving of ${}^pN_0^V$.
This data indicates a dependence of accuracy on ${}^pN_0^V$, which was not recognized in the analysis of Brandt and Diskin.
These are the first published multilevel numerical results of a particular segmental refinement algorithm; new algorithms and implementations should be developed with, for instance, higher order prolongation and vertex-centered discretizations to determine if this dependence can be ameliorated or if it is a fundamental property of the method.

\section{Segmental refinement communication complexity}
\label{sec:model}

Segmental refinement inherits the computational depth of conventional distributed memory multigrid, and the coarse grid solves are identical; a more refined complexity model is required to distinguish the communication characteristics of SR from those of conventional multigrid.
Mohr has analyzed the communication patterns and savings with SR and the extra computation costs for a two-level SR method \cite{Mohr1998,Mohr2000}.
Brandt has presented memory complexity analysis with a $\log^D$ term for the SR buffer cells memory complexity (\S 8.7 \cite{ABrandt_1984b}).
This section proposes a new SR data model, that we posit is asymptotically exact, and an abstract memory model that resolves the communication that is eliminated by this new SR method.

% to better understand the theoretical implications of segmental refinement on multigrid communication complexity.

\subsection{A multigrid V--cycle communication model}
\label{ssec:comm_model_v}

We define two types of multigrid communication: vertical inter-grid ($c_V$) and horizontal intra-grid ($c_H$) communication.
%These two types of communication have different character in both conventional and SR multigrid; vertical communication is volumetric data and horizontal communication is surface data.
%Horizontal communication presents the opportunity to buffer and aggregate messages in, for instance, the smoother and residual computation.
%Message aggregation also allows for loop fusion to increase the arithmetic intensity of conventional  multigrid components \cite{Williams2012}.
%Segmental refinement can be viewed as applying message aggregation to the entire set of finest grids, but without replicating the semantics exactly.
A conventional distributed memory multigrid V--cycle uses 26 ($c_H$) messages per process in each residual, smoother, and operator application when using a 27-point stencil in 3D with a Cartesian process grid in a standard nearest-neighbor exchange process.
A $V(2,2)$ cycle requires 156 horizontal messages per level, including the $\tau$ correction term, plus vertical restriction and prolongation messages in eight message phases, or bulk synchronous steps (six horizontal and two vertical).
Our model focuses on these communication phases and the ``distance" of each message phase.

\subsection{A memory model} 
\label{sec:machine}

We define a ``word" of data as a small patch of cells (e.g.,  $4^D-32^D$ cells) and assume that each ``process" computes on one data word.
Consider a two-level memory model with $Q$ words of fine grid memory, partitioned into $\sqrt{Q}$ partitions, each of size $\sqrt{Q}$.
Level $1$ is on-partition memory and level $2$ is off-partition memory. 
We define near communication as communication between processes within a memory partition and far communication as communication between memory partitions.

\subsection{Proposed asymptotic segmental refinement data model}
\label{sec:new_sr_model}

%\begin{figure}[h]
\begin{wrapfigure}{R}{0.5\textwidth}
\begin{center}
\includegraphics[width=80mm]{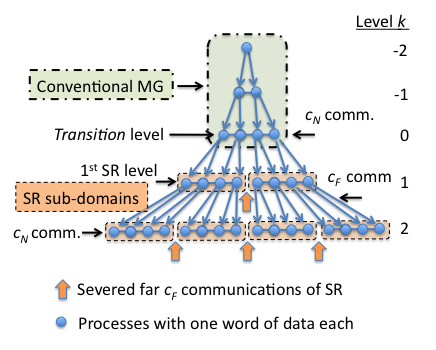} 
\caption{1D example of asymptotic data model}
\label{fig:1dsr}
\end{center}
%\end{figure}
\end{wrapfigure}
%Figure \ref{fig:1dsr} (a) shows the SR data model as implemented and tested in \S \ref{ssec:exp_params} and \S\ref{sec:numr}, and Figure \ref{fig:1dsr} (b) shows the same problem with the data model analyzed in this section.

%A simple and practical use of SR would be to fix $K=5$; this reduces the size of the full communication (coarse grid) solver by a factor of 32K with respect to the traditional parallel multigrid method.
%However, we focus on asymptotics for analysis.
The observations in \S\ref{ssec:exp_params} suggest that a data model that increases ${}^pN_0^V$ with the number of levels, and perhaps adds a quadratic term to equation (\ref{eq:nb}), would be asymptotically exact.
One could attempt to keep ${}^pN_0^V$ constant and determine an appropriate buffer schedule, but this would be less reliable and natural given our understanding of the problem.
Note, a non-asymptotic model could be useful in practice because, for instance, fixing $K=5$ reduced the size of the conventional (full communication) solver by a factor of 32K ($2^{KD}$), which is a significant constant.

We propose a data model that we posit would provide sufficient accuracy for an asymptotically exact solver, by extending the parallel octree of the coarse grids to the entire multigrid hierarchy, using $\sqrt{Q}$ processes in each SR subdomain and setting the size of the transition level to fit into one memory partition.
With $K$ SR levels, this model has $M=2K$ multigrid levels and $K+1$ conventional levels. 
%Refine this single subdomain of the transition level to create eight ($2^D$) SR patches of the same logical size, populating eight memory partitions, and continue this process until all partitions are populated with $L=2K+1$ multigrid levels, $K$ SR levels, and $K+1$ conventional levels. 
Figure \ref{fig:1dsr} shows a 1D example of this data model with two SR levels, $Q=16$, and an SR patch length $N_0^V=4$ where a word is one cell.

\subsection{Communication complexity}
\label{ssec:comm}

We use the multigrid V--cycle communication model of \S\ref{ssec:comm_model_v}, the machine model of \S\ref{sec:machine}, and the SR data model of \S\ref{sec:new_sr_model} to analyze the communication complexity this segmental refinement method.
We ignore FMG prolongation because there are $M$ FMG prolongations as opposed to \Order{M^2} V--cycle restrictions and prolongations.
FMG processes a V--cycle once on the finest grid, twice on the first coarse grid with $1/2^{D}$ as many active processes, and so on for $M+1$ levels and $\left(M+1\right)\cdot\frac{M}{2} \approx \frac{M^2}{2}$ grid visits total.
This is the source of the computation depth, $log^2_2 \left(N\right)$, of FMG.
There are six $c_H$ communication phases and two $c_V$ phases per grid visit with one visit on the finest grid, two on the first coarse grid, and so on with $M$ visits to the coarsest grid.
There are about $M^2/8=M^2/2^D$ visits on the fine (SR) half of the grid hierarchy and $3\cdot M^2/8$ on the coarse (conventional) half of the grid hierarchy.
We ignore vertical data locality and assume that all vertical communication is far communication in the finest $K$ levels and that all communication is near communication on the coarsest $K+1$ levels. 
%All horizontal communication is near communication in the SR solver;  the conventional solver requires far communication on the finest $K$ levels.

\subsubsection{3D Bisection bandwidth}
\label{ssec:bbw}

Briefly consider a four level memory model generated by bisecting the memory and domain of the current model.
The communication complexity between these two partitions is bisection bandwidth.
The highest order term of bisection bandwidth of conventional multigrid is from the ghost cell exchange on the finest grid.
On a 3D cube with $N$ cells in each dimension, \Order{N^2} cells (the area of the face between the two partitions) times the length of the ghost region (\Order{1}) is communicated \Order{1} times resulting an a communication complexity of \Order{N^2}.

The highest order term in the SR bisection bandwidth complexity is from the buffer region exchange on the transition level.
% if we assume that a local buffer for the domain of $I_0^{1}$ is filled with a horizontal exchange on level $k=0$.
%the finest grid that communicates between two halves of the machine and is between two sets of four $\sqrt Q$ partitions each. 
Assume the number of buffer cells required is quadratic in $K$, because our data in \S\ref{ssec:exp_params} suggests this is required for a convergent solver.
The ``area" of data sent in this buffer cell exchange is $\Order{{}^pN_0^2}=\Order{\sqrt{N}^2}=\Order{N}$; it has a depth $K^2$ and is executed \Order{\log_2{N}} times.
Thus, the communication complexity is $ \Order{N \cdot K^2}\Order{\log_2{N}} =  \Order{N \cdot \log_2^3{N}}$.
Segmental refinement reduces the bisection communication requirements from $\Order{N^2}$ to $\Order{N\log_2^3{N}}$.

\subsubsection{Near and far communication complexity}
\label{ssec:commcomplex}

Table \ref{tab:comm} tabulates the communication complexity of conventional and SR multigrid with $M+1$ levels. %, $K$ SR levels ($L=2K+1$).
The coarsest $K+1$ levels of both solvers use the same FMG solver on one memory partition.
There are six $c_H$ communication phases and two $c_V$ phases per grid visit. %;  $L^2/8$ visits on the fine half of the grid hierarchy and $3\cdot L^2/8$ visits on the coarse half.
%Segmental refinement removed the $6c_H$ term from the far communication complexity -- this is the distinguishing characteristic of segmental refinement.
\begin{table}[h]
\begin{center}
\begin{tabular}{|c||c|c|}
\hline
Communication type   &  Near &  Far \\ \hline
Coarse grids  &  $3\cdot\left(6c_H + 2c_V\right) $          &  $0$  \\
Conventional fine grids&  $6c_H$          &  ${6c_H} + 2c_V$   \\
SR fine grids  &  $6c_H$          &  $2c_V$   \\
\hline
\end{tabular}
\caption{Communication phases ($\times \log^2_2{N}/8$) of conventional distributed memory multigrid and segmental refinement multigrid}
\label{tab:comm}
\end{center}
\end{table}
The removal of far horizontal communication complexity is the distinguishing characteristic of segmental refinement.

A given segmental refinement data model removes horizontal communication at some level of the memory hierarchy, the $6c_H$ term in far communication in Table \ref{tab:comm} and at the arrows Figures \ref{fig:1dsr-a} and \ref{fig:1dsr}.
Communication, in some memory model, is used only for the vertical operators restriction and prolongation, which have tree-like graphs.
Tree algorithms are efficient for the global communication required for the solve of an elliptic system.
The critical observation of segmental refinement is that horizontal communication of traditional parallel multigrid is used for local processes and is not global, hence ``far" communication is potentially not necessary.
%Removing horizontal communication, at some scale, is the distinguishing characteristic of segmental refinement with respect to communication complexity.

We have investigated two segmental refinement data models that are ``two level" in that there is one ``transition" level between a conventional coarse grid solver and decoupled finer grids.
One can, in principle, compose these two models, by using the method in this section as the coarse grid solver for the method in \S\ref{sec:sr}, and create a three level method.
We speculate that one could generate an asymptotically exact ``multilevel" segmental refinement method that starts with the method in this section as a ``coarse grid" solver and reduces the size of ${}^pN^V_k$, by a factor of two on each finer level, resulting in just one process per SR subdomain on a fine level, thereby recovering more parallelism, and continue with the method in \S\ref{sec:sr} on each process.
This is a subject for future work.

\section{Timing studies}
\label{sec:numr}

This section presents scalability data on the problem in \S\ref{ssec:problem} on the Cray XC30 at NERSC, with up to 64K cores.
We use  8 of the 12 cores on each socket and thus utilize 96K cores at scale, or about 75\% of the machine and investigate weak scaling with $128^3$ and $32^3$ cells per core on the fine grid, with four and three SR levels respectively (and ${}^pN_0^V=8$ and $4$ respectively).
The solver is preloaded with one solve, which verifies accuracy, followed by 8 timed solves for the $128^3$ cells per core case and 512 solves for the $32^3$ cells per core case to normalize times.

\begin{figure}[h]
\begin{center}
\begin{minipage}[t]{.48\linewidth}
\includegraphics[width=88mm]{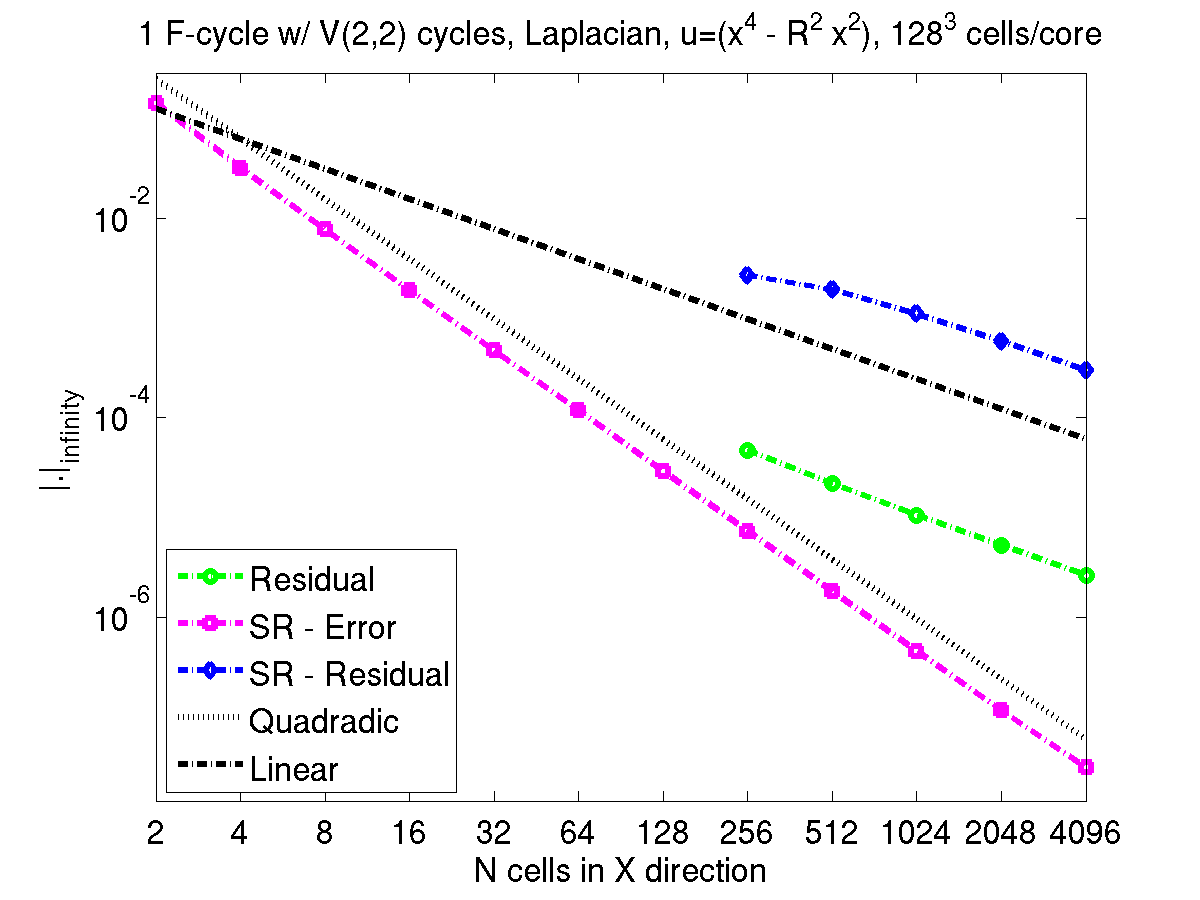} 
\caption{Convergence verification}
\label{fig:conva}
\end{minipage}
\begin{minipage}[t]{.48\linewidth}
\includegraphics[width=88mm]{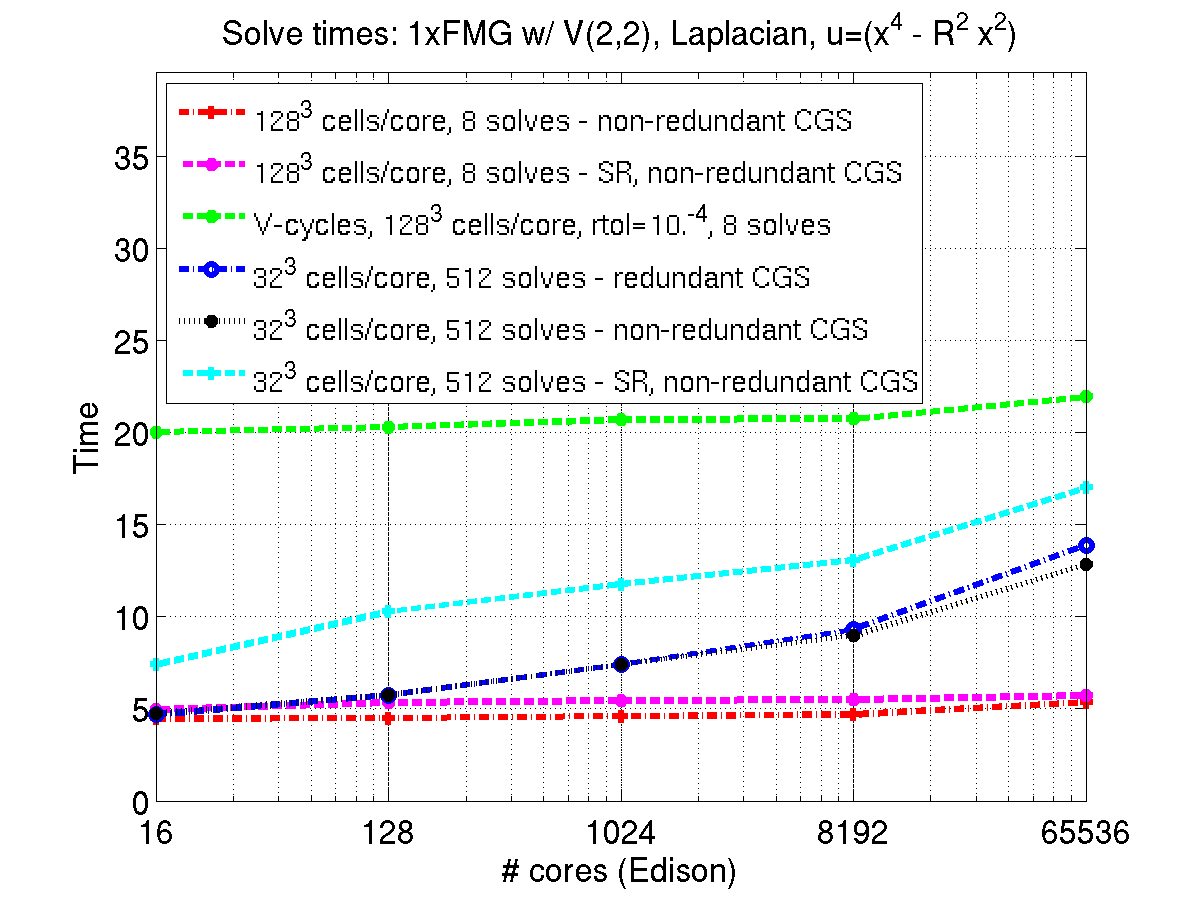} 
\caption{Edison weak scaling}
\label{fig:scaling}
\end{minipage}
\end{center}
\end{figure}

\begin{wrapfigure}{R}{0.5\textwidth}
%\begin{figure}[ht]
\begin{center}
\includegraphics[width=88mm]{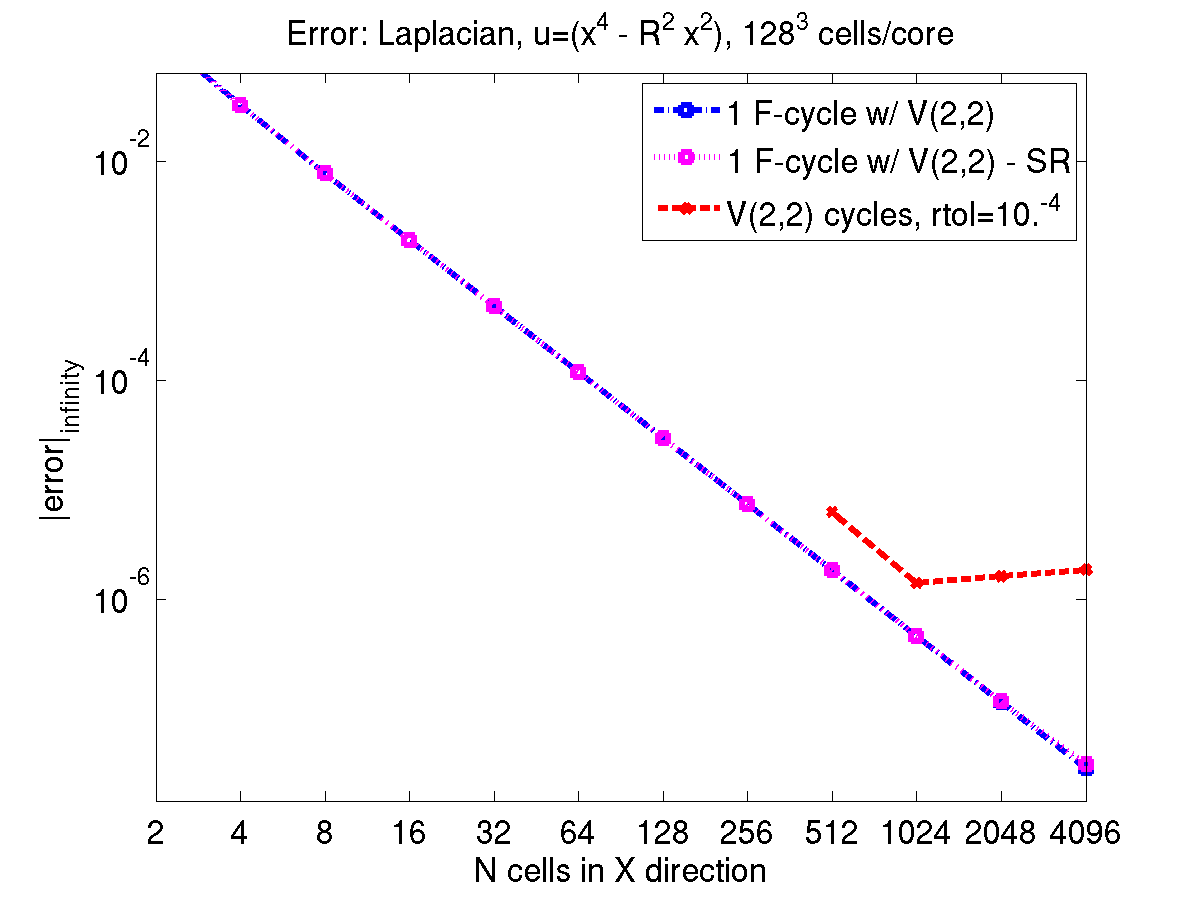} 
%\captionsetup{justification=raggedright,singlelinecheck=false}
\caption{Errors of all solvers}
\label{fig:verrors}
\end{center}
%\end{figure}
\end{wrapfigure}

Figure \ref{fig:conva} plots the infinity norm of the error and residual in the FMG solve and verifies that our solvers are asymptotically exact and that $2^{nd}$-order accuracy is achieved, but only $1^{st}$-order reduction is observed in the residual.
The residuals for SR are larger than those of the conventional method but are still $1^{st}$-order convergent.
Figure \ref{fig:scaling} plots the solve times for the SR solver and the conventional multigrid solvers and shows modest gains in scalability with SR.
The solve times for a V--cycle solve with a relative residual tolerance of $10^{-4}$ are also shown.
Figure \ref{fig:verrors} demonstrates the stagnation in error reduction with a V--cycle solver, converged to a constant residual reduction, and that SR is maintaining perfect $2^{nd}$-order accuracy.

Coarse grids can be computed redundantly, where all processors are active on all levels redundantly computing coarse grid corrections, or processors can be left idle on coarse grids.
Redundant coarse grid solves result in a ``butterfly" communication pattern and the idle processors result in a tree communication pattern.
This approach has the advantage of requiring no communication in the prolongation phase, hence reducing the number of bulk synchronous communication steps, at the expense of sending more data with more messages overall.
We observe in Figure \ref{fig:scaling} that redundant coarse grid solves are slightly slower, which suggests that larger number of messages cost more than the savings in the number of bulk synchronous phases.
This could be due to contention in the network during restriction, however, the differences are small and only observable on the largest run.

%\clearpage
\section{Conclusions}
\label{sec:conc}

We continue the work of Brandt and Diskin \cite{ABrandt_BDiskin_1994a}, with the first published multilevel numerical results of the segmental refinement multigrid method.
We demonstrate that SR can maintain the semantics of textbook efficient multigrid FMG-FAS with processing that is more attractive on modern memory-centric architectures than conventional distributed memory multigrid by decoupling fine grid processing, which improves data locality; amortizes latency costs; and reduces data dependencies. 
We have experimentally investigated the asymptotic behavior of SR and have found an accuracy dependance not previously recognized.
We analyzed the communication complexity, with a two-level memory model, of an SR data model, where we show that the method removes horizontal communication as define by the memory model.
The degree to which the memory model, on which any given SR method removes communication, is a useful performance model for any given machine can be used as a metric for the potential efficacy of the method.
We experimentally verify that our SR data model is an asymptotically exact solver on 64K cores of a Cray XC30 and provide timing and scaling data.

We have observed modest improvement in scaling with SR with a simple data model that supports only a few SR levels.
Future work includes developing SR data models that accommodate more levels of the memory hierarchy, testing on machines with deeper memory hierarchies and fully exploiting SR's data locality with, for instance, loop fusion \cite{Williams2012}.
A vertex-centered discretization and high order $I_0^1$ prolongation would be of interest to better understand the asymptotic complexity of SR and corroborate the observation of the accuracy dependance on ${}^pN^V_0$.
SR may be particularly sensitive to the order of prolongation because it is used to set the ``frozen" ghost cells in the SR buffer region.
We have investigated a model problem; further work involves extending the application of SR to more domains, such as variable coefficient and nonlinear problems and unstructured grid problems.

All code, data, and run and parsing scripts used in this paper are publicly available at https://bitbucket.org/madams/srgmg.

\subsection*{Acknowledgments}
We would like to thank Achi Brandt for his generous guidance in developing these algorithms.
This material is based upon work supported by the U.S. Department of Energy, Office of Science, Office of Advanced Scientific Computing Research and performed under the auspices of the U.S. Department of Energy by Lawrence Berkeley National Laboratory under Contract DE-AC02-05CH11231.
This research used resources of the National Energy Research Scientific Computing Center, which is a DOE Office of Science User Facility.
Authors from Lawrence Berkeley National Laboratory were supported by the U.S. Department of Energy's Advanced Scientific Computing Research Program under contract DEAC02-05CH11231.

%\appendix
%\section{How I became inspired}

%(RS: reference TBD01 -- check ``navierstokes'' and ``reynolds'' -- you can capitalize in the bib file by putting them in curly braces i.e. {Reynolds}.) MA - done

\bibliographystyle{amsalpha}
\bibliography{./bib}

\newcommand{\etalchar}[1]{$^{#1}$}
\providecommand{\bysame}{\leavevmode\hbox to3em{\hrulefill}\thinspace}
\providecommand{\MR}{\relax\ifhmode\unskip\space\fi MR }
% \MRhref is called by the amsart/book/proc definition of \MR.
\providecommand{\MRhref}[2]{%
  \href{http://www.ams.org/mathscinet-getitem?mr=#1}{#2}
}
\providecommand{\href}[2]{#2}
\begin{thebibliography}{WKS{\etalchar{+}}12}

\bibitem[ABS{\etalchar{+}}14]{HPGMGv1}
M.~F. Adams, J.~Brown, J.~Shalf, B.~Van Straalen, E.~Strohmaier, and
  S.~Williams, \emph{{HPGMG} 1.0: A benchmark for ranking high performance
  computing systems}, Tech. Report LBNL-6630E, LBNL, Berkeley
  (https://bitbucket.org/hpgmg/hpgmg), 2014.

\bibitem[ASB10]{Adams-10a}
M.~F. Adams, R.~Samtaney, and A.~Brandt, \emph{Toward textbook multigrid
  efficiency for fully implicit resistive magnetohydrodynamics}, Journal of
  Computational Physics \textbf{229} (2010), no.~18, 6208 -- 6219.

\bibitem[BD81]{Bank-81}
R.E. Bank and T.~Dupont, \emph{An optimal order process for solving finite
  element equations}, Math. Comp. \textbf{36} (1981), 35--51.

\bibitem[BD94]{ABrandt_BDiskin_1994a}
A.~Brandt and B.~Diskin, \emph{Multigrid solvers on decomposed domains}, Domain
  Decomposition Methods in Science and Engineering: The Sixth International
  Conference on Domain Decomposition (Providence, Rhode Island), Contemporary
  Mathematics, vol. 157, American Mathematical Society, 1994, pp.~135--155.

\bibitem[BL11]{Brandt-2011}
A.~Brandt and O.~E. Livne, \emph{Multigrid techniques}, Society for Industrial
  and Applied Mathematics, 2011.

\bibitem[Bra73]{ABrandt_1973a}
A.~Brandt, \emph{Multi--level adaptive technique ({MLAT}) for fast numerical
  solution to boundary value problems}, Proceedings of the Third International
  Conference on Numerical Methods in Fluid Mechanics (Berlin) (H.~Cabannes and
  R.~Teman, eds.), Lecture Notes in Physics, vol.~18, Springer--Verlag, 1973,
  pp.~82--89.

\bibitem[Bra77]{Brandt-77}
A.~Brandt, \emph{Multi-level adaptive solutions to boundary value problems},
  Math. Comput. \textbf{31} (1977), 333--390.

\bibitem[Bra84]{ABrandt_1984b}
A.~Brandt, \emph{Multigrid techniques: 1984 guide with applications to fluid
  dynamics}, GMD--Studien Nr. 85, Gesellschaft f{\"u}r {M}athematik und
  {D}atenverarbeitung, St. Augustin, 1984.

\bibitem[Din79]{NDinar_1979a}
N.~Dinar, \emph{Fast methods for the numerical solution of boundary value
  problems}, Ph.D. thesis, Weizmann Institute of Science, Rehovot, Isreal,
  1979.

\bibitem[Fed61]{RPFedorenko_1961a}
R.~P. Fedorenko, \emph{A relaxation method for solving elliptic difference
  equations}, Z. Vycisl. Mat. i. Mat. Fiz. \textbf{1} (1961), 922--927, Also in
  U.S.S.R. Comput. Math. and Math. Phys., 1 (1962), pp. 1092--1096.

\bibitem[Moh00]{Mohr2000}
M.~Mohr, \emph{Low communication parallel multigrid}, Euro-Par 2000 Parallel
  Processing (Arndt Bode, Thomas Ludwig, Wolfgang Karl, and Roland Wismüller,
  eds.), Lecture Notes in Computer Science, vol. 1900, Springer Berlin
  Heidelberg, 2000, pp.~806--814 (English).

\bibitem[MR98]{Mohr1998}
M.~Mohr and U.~Rude, \emph{Communication reduced parallel multigrid: Analysis
  and experiments}, Tech. Report Technical Report No. 394, University of
  Augsburg, 1998.

\bibitem[Sou40]{RVSouthwell_1940a}
R.~V. Southwell, \emph{Relaxation methods in engineering science}, Oxford
  University Press, Oxford, 1940.

\bibitem[TDB01]{Thomas2001853}
J.~L. Thomas, B.~Diskin, and A.~Brandt, \emph{Textbook multigrid efficiency for
  the incompressible {N}avier--{S}tokes equations: high {R}eynolds number wakes
  and boundary layers}, Computers \& Fluids \textbf{30} (2001), no.~78, 853 --
  874.

\bibitem[TOS01]{UTrottenberg_CWOosterlee_ASchueller_2000a}
U.~Trottenberg, C.~W. Oosterlee, and A.~Sch{\"u}ller, \emph{Multigrid},
  Academic Press, London, 2001.

\bibitem[WKS{\etalchar{+}}12]{Williams2012}
S.~Williams, D.~D. Kalamkar, A.~Singh, A.~M. Deshpande, B.~Van~Straalen,
  M.~Smelyanskiy, A.~Almgren, P.~Dubey, J.~Shalf, and L.~Oliker,
  \emph{Optimization of geometric multigrid for emerging multi- and manycore
  processors}, Proceedings of the International Conference on High Performance
  Computing, Networking, Storage and Analysis (Los Alamitos, CA, USA), SC '12,
  IEEE Computer Society Press, 2012, pp.~96:1--96:11.

\end{thebibliography}

\end{document}